\documentclass[a4paper]{article}

\usepackage{amsfonts,amssymb}
\usepackage{a4wide}
\usepackage{graphicx,geompsfi-pdf}
\usepackage{subfig}
\usepackage{caption}
\usepackage{multirow}
\usepackage{afterpage}

\newtheorem{theorem}{Theorem}

\def\R{\mathbb{R}}

\begin{document}

\title{Applied Mathematics, the \emph{Hans van Duijn} way}
\author{Mark A. Peletier}
\date{\today}
\maketitle

\begin{abstract}
This is a former PhD student's take on his teacher's scientific philosophy. I describe a set of `principles' that I believe are conducive to good applied mathematics, and that I have learnt myself from observing Hans van Duijn in action. 
\end{abstract}

\section{Introduction}

Among scientists there circulates an apocryphal list of six `principles of good science', supposedly formulated by Clifford Truesdell. My copy is a bad photocopy of a bad photocopy of this list. It is entitled `How the Bernoullis worked', and, despite its obscure provenance, gives some very sound advice: `Always attack a specific problem. Try to solve it in a way that leads to a general method', for instance, and `Let a key problem solved be father to a key problem posed'. 

This list appears inspired by the belief that there are `good' and `bad' ways of doing science. Over the years I have come to agree. I also claim that similarly there are good and bad ways of doing \emph{applied mathematics}, and that the peculiarities of applied mathematics warrant a special treatment. 

The list below is strongly inspired by my experience of working with Hans van Duijn, and in this volume dedicated to Hans' work I am perfectly confident giving it the title `How Hans works'. The list is wholly my own doing, however, and Hans should not be blamed for any (undoubtably many) faults. In fact, the occasions on which he and I discussed such philosophical issues were few and far between. This list is my analysis of Hans' actions, rather than his words---a written version of `leading by example'.

\bigskip

{\centering
\fboxsep15pt
\fbox{\sffamily\begin{minipage}{0.7\textwidth}
\centering {\large How Hans works}
\begin{enumerate}
\item Find a math-friendly applier
\item Get to understand the applied problem \emph{completely}
\item In the applied mess, find a beautiful math problem 
\item Solve it, preferably \emph{with} the appliers 
\item Interpret it, \emph{with} the appliers 
\item Repeat and enjoy!
\end{enumerate}
\end{minipage}}

}

\bigskip

In the rest of this paper I will illustrate these `rules' by describing some scientific developments that I have been part of.

\section{Mathematical Geology}

The topic is what one might call `Mathematical Geology'. Central questions are for instance
\begin{itemize}
\item Forward problems: how do rocks deform? How is this deformation influenced by material properties and loading conditions?
\item Inverse problems: from observations of current rocks, what can we infer  about the deformation history?
\end{itemize}
The inverse problem has of course significant commercial interest, although our mathematical work never seemed to get close to those interests.

\begin{figure}[t]
\centering
\includegraphics[height=7cm]{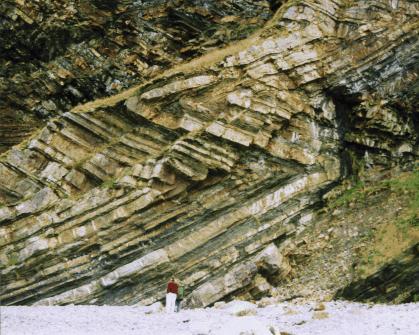}
\caption{An outcrop at Millook Haven in Cornwall, UK. The lighter layers are sandstone, the darker ones shale (clay). At the bottom of the picture there is a young version of myself.}
\label{fig:MillookHaven}
\end{figure}

\begin{figure}[ht]
\centering
\includegraphics[height=6cm]{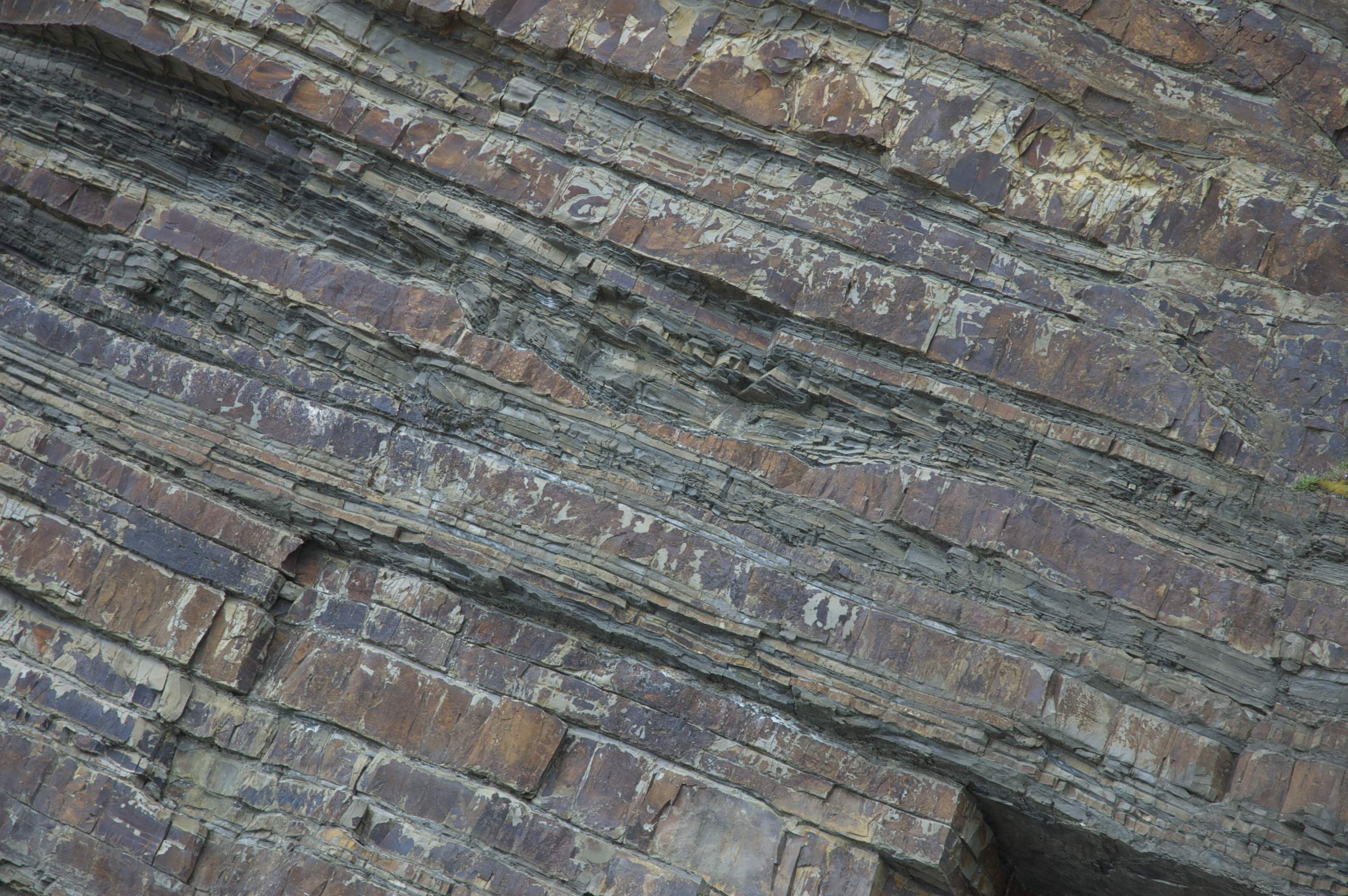}
\caption{Sandstone and shale, similar to Fig.~\ref{fig:MillookHaven}. Here the sandstone has rusted and become brown. The scale from top to bottom is about two meters. This photograph was taken to highlight the unusual \emph{deviation} from straight, parallel layers: exactly in the middle of the picture a sandstone layer has a thinner section, formed by a small channel of moving water eroding a layer of sand just after deposition.}
\label{fig:channel}
\end{figure}

The geological structures that we study were mostly created by sedimentation, which creates parallel layers of rock whose thickness and composition vary from one layer to the next but are  remarkably constant along each layer (see Fig.~\ref{fig:channel}). By the time we see these layers, as in Fig.~\ref{fig:MillookHaven}, they have gone through a long series of deformations: submersion to large depth, conversion of the sediment into rock, folding under various forces, possibly in many stages, and possibly combined with various chemical modification processes as well. The whole deformation history can be very complex, and we tend to focus on just one part: the folding of the layers.

Much of what I know about mechanics and rock folding I have learnt from other scientists, which brings me to Rule \#1: \emph{Find a math-friendly applier}. The importance of this rule can not be overemphasized. There have been many such applied collaborators over the years, but here I want to illustrate this principle with two of these: Giles Hunt and John Cosgrove. 

\medskip
\begin{center}
\noindent\vbox{\halign{\hfil #\hfil& \qquad \hfil # \hfil\cr\includegraphics[height=5cm]{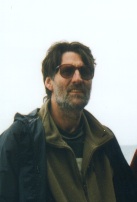}&\includegraphics[height=5cm]{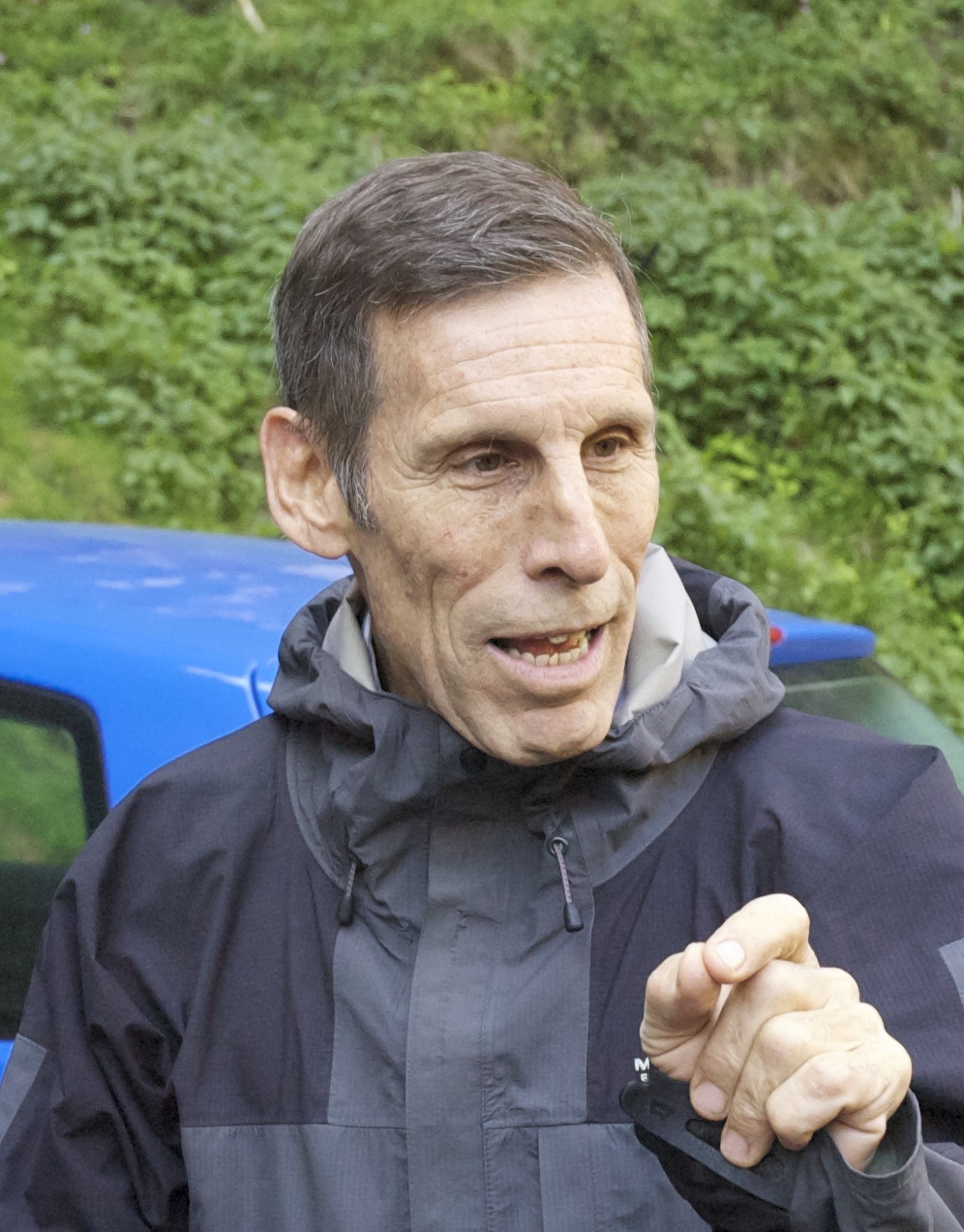}\cr Giles Hunt & John Cosgrove\cr}}
\end{center}

Giles Hunt is a retired mechanical engineer at the University of Bath, and I was fortunate enough to spend a postdoc period in his department. John Cosgrove is a geologist at Imperial College in London whom I met at various occasions, and who brought the rocks to life for us on a field trip in the spring of 2011.

What makes both Giles and John such great collaborators for a mathematician is a combination of three things. First: an untiring willingness to explain how they think about `their' problems. Anyone who has ever worked outside of their field of expertise will appreciate the importance of this. Secondly, the ability to think abstractly, to seek general principles that govern the specific case at hand. Finally, the belief that mathematics, and more generally reductionism, might help them understand their problems better. Oh, and of course enthusiasm always helps. 

To give an example, I spent many hours with Giles experimenting with toy models of mechanical systems, which were contraptions of links, hinges, and springs---meccano on paper. Giles is a master in using such highly simplified mechanical models to zoom in on a very specific feature while preserving mechanical feasibility. All the models we built on paper could, in theory, be built in a lab. Some actually were. 

\medskip

\begin{figure}[p]
\includegraphics[width=\textwidth]{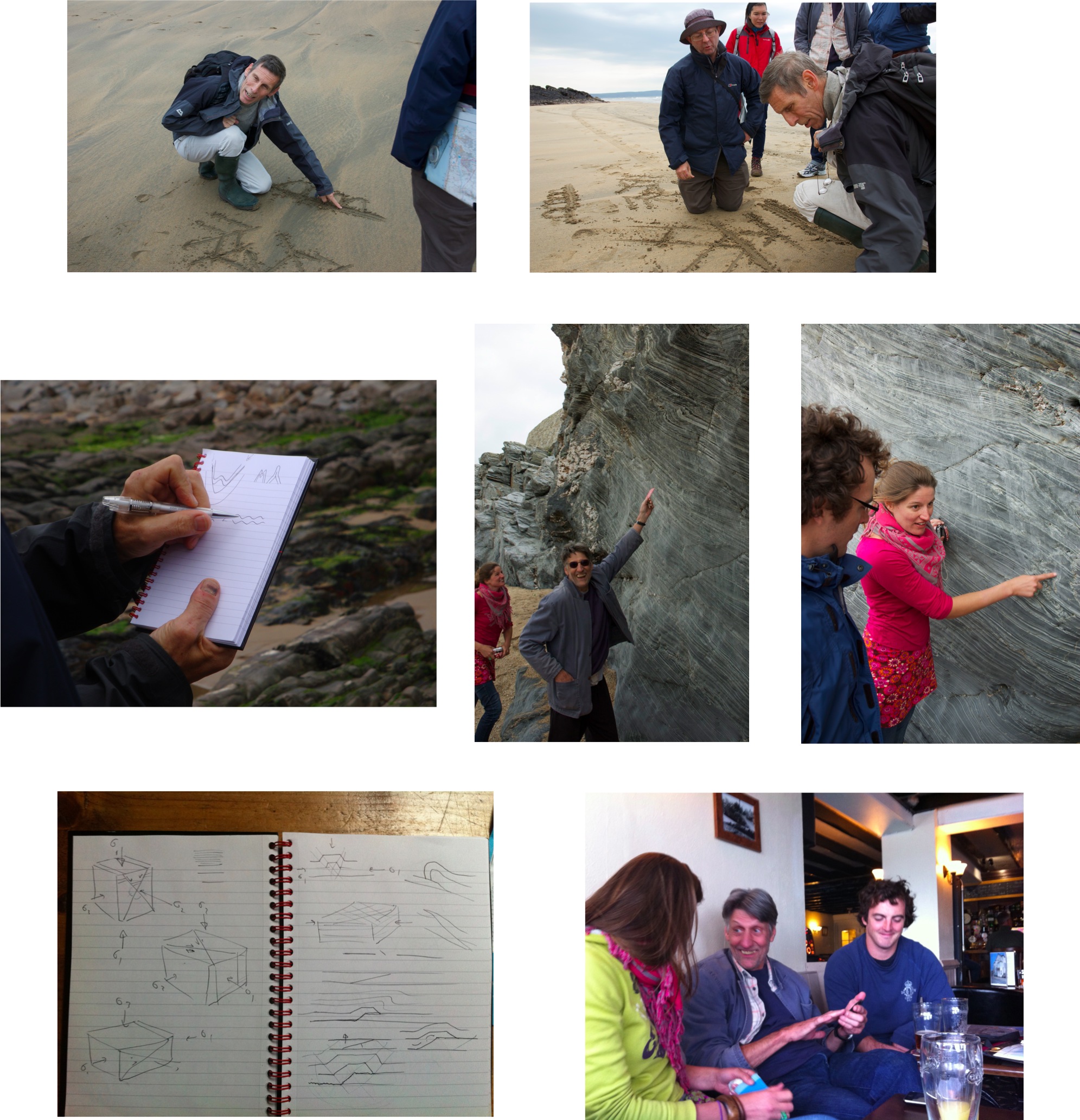}
\caption{An impression of a field trip to Cornwall}
\label{fig:field-trip}
\end{figure}

This brings me to Rule \#2: \emph{Get to understand the problem completely}. In the case of mathematical geology, there is no substitute to a field trip with a guide like John Cosgrove. Over the course of a couple of days we saw many different types of rocks with very varying deformation histories, and it was John's expert insight that made them accessible to us. Figure~\ref{fig:field-trip} shows some of the activities during the trip. 

In a similar way, an earlier field trip, in 1998, had pointed us towards the first of a number of mathematical questions: 

\begin{quote}
\textbf{Math question \#1: } Why straight limbs and sharp hinges?
\end{quote}
Why this is interesting requires a bit of background. Rocks are fairly rigid at human time scales, but over geological time scales they can deform significantly, as the deformed rocks around the world show. Whether rock deforms elastically or viscously on geological time scales is still unclear, to large extent due to the difficulty to performing experiments at those time scales. But both elastic and viscous behaviour have a preference for smooth deformation that is spread out, not localized. The deformation in Fig.~\ref{fig:MillookHaven} however is highly localized: long stretches of undeformed rock are bounded by very sharp corners. Why is this?

\subsection{The first mathematical problem}

Math question \#1 is a very intriguing one, and one of the first to come out of our interest in geology (see Rule \#3: \emph{In the applied mess, find a beautiful math problem}). It is also a tough one, and not yet completely resolved. 

\begin{figure}[htb]
\centering
\includegraphics[height=4cm]{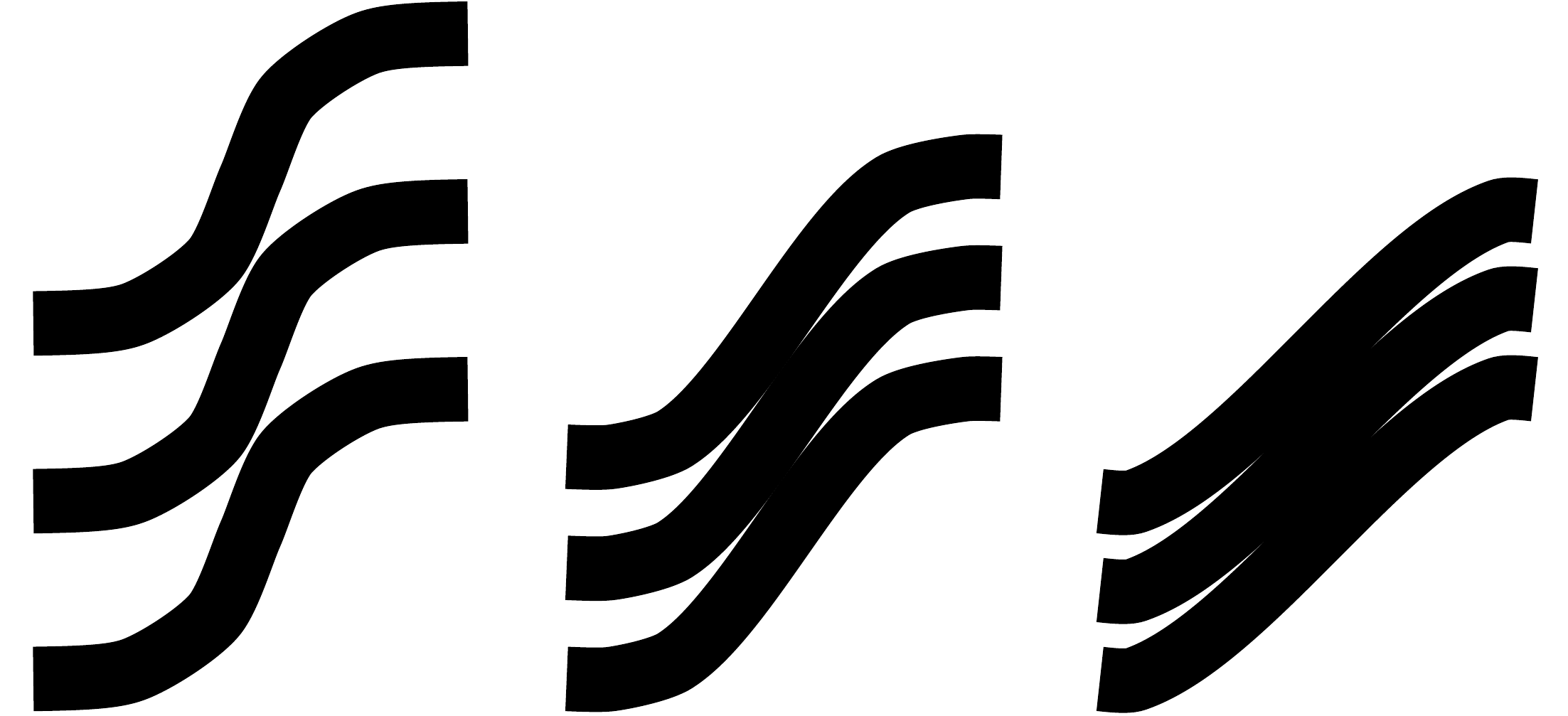}
\caption{A pictorial answer to Math question \#1: high overburden pressure implies that voids between the layers are penalized, leading to straight limbs with sharp corners. Published in a paper with Giles Hunt and Ahmer Wadee~\cite{HuntPeletierWadee00} (see Rule \#4: \emph{Solve it, preferably with the appliers}).}
\label{fig:sharp-corners}
\end{figure}

After that first field trip we did give an answer of sorts. Not so much a mathematical answer, but more a pictorial one. By simply fitting together curved layers of different geometry (Fig.~\ref{fig:sharp-corners}) we `showed' that the straight-limb-sharp-corner geometry fits together better, leaving less voids between the layers. Since the layers deform under high pressure, voids carry a high energy penalty. This penalty pushes the material away from its preference for smooth deformation towards the sharp-hinged geometry.

\begin{figure}[htb]
\centering
\includegraphics[width=0.9\textwidth]{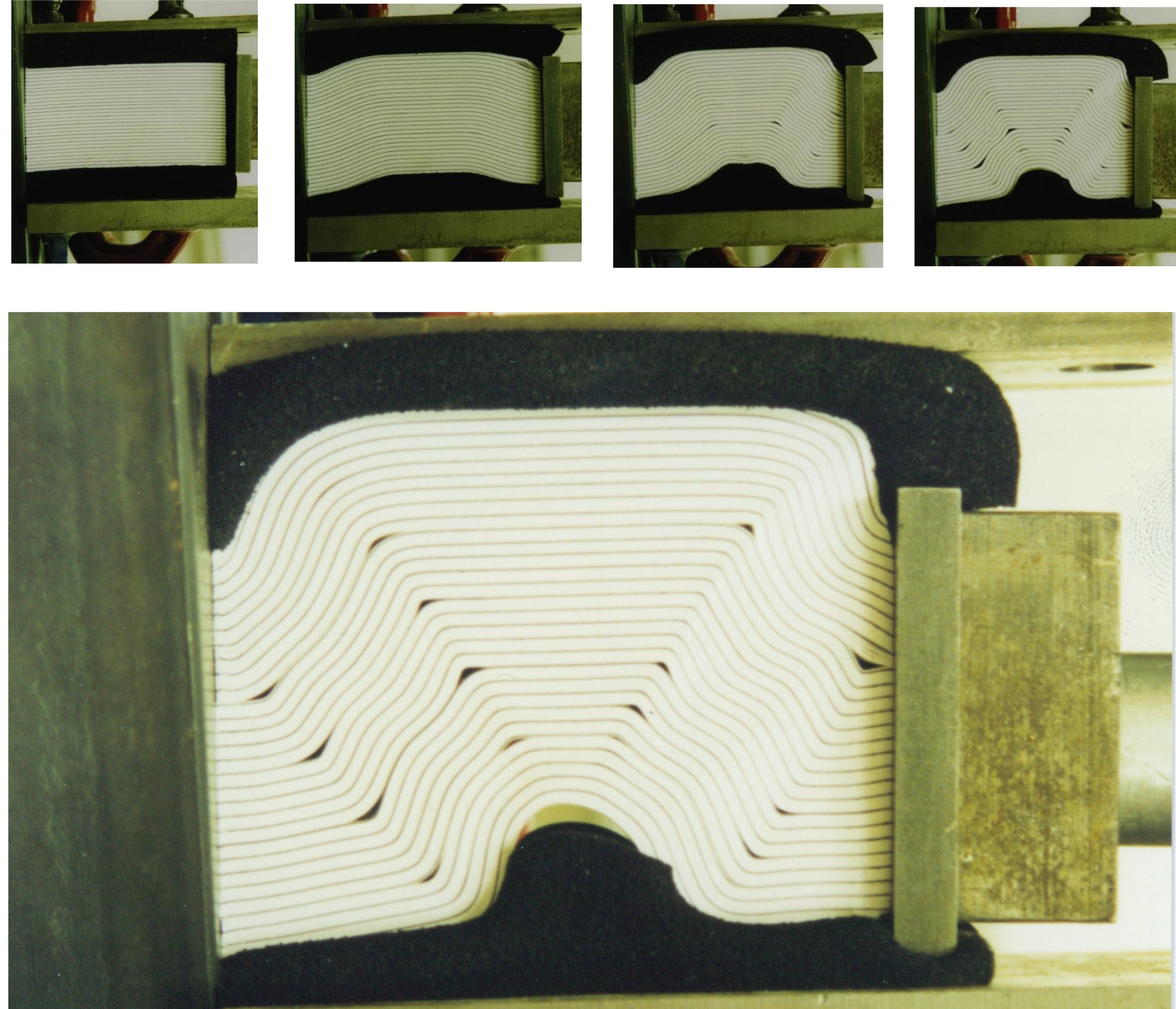}
\caption{Horizontal compression of layers of paper, vertically confined by two layers of foam. The top four figures show the compression increasing from left to right; the bottom picture is the final structure. One layer in twenty-five has been coloured black. These experiments were done with Giles Hunt and Ahmer Wadee at the University of Bath.}
\label{fig:c5}
\end{figure}

Later we would come back to this question. In the meantime Bruce Hobbs, another geologist, had suggested that we do some experiments in the lab, using paper as a substitute for layered rocks. Figure~\ref{fig:c5} shows an example of this. 

\afterpage{\clearpage}

The setup of the experiment is as follows. We placed a stack of (relatively common, printer-type) paper between two confining plates, lined with foam. We then slowly compressed the layers in the direction of the layering. As as result, the stack of paper buckled into the foam; as the shortening increased, the foam became stiffer, since it was flanked by rigid plates, and the smooth initial buckle converted into the final  structure in Fig.~\ref{fig:c5}.

\begin{figure}[htb]
\centering
\includegraphics[height=3cm]{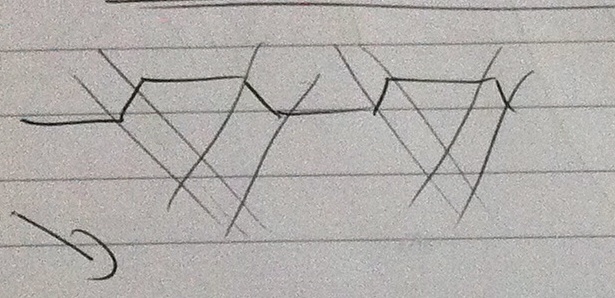}
\qquad
\includegraphics[height=3.5cm]{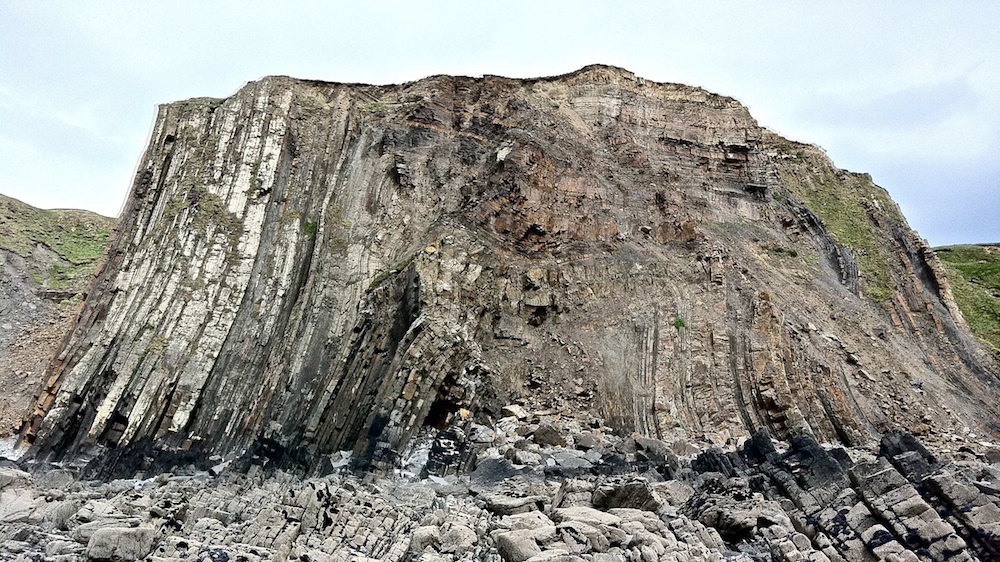}
\caption{The left-hand figure shows John Cosgrove's schematic of `conjugate kink bands', and the right-hand figure an example of such structures in Bude.}
\label{fig:kinkbands}
\end{figure}

This experiment showed us a number of things. First, there is indeed a tendency for the paper to be straight with sharp corners, and this tendency increases with increasing shortening. Since the horizontal shortening also causes the foam to compress, the total pressure in the paper increases, which explains why the sharpness of the corners increases with shortening. It also showed us a structure that appears often in geology, that of `kink bands' (see Fig.~\ref{fig:kinkbands}). 

\subsection{Kink bands}

\begin{figure}[htb]
\centering
\includegraphics[height=6cm]{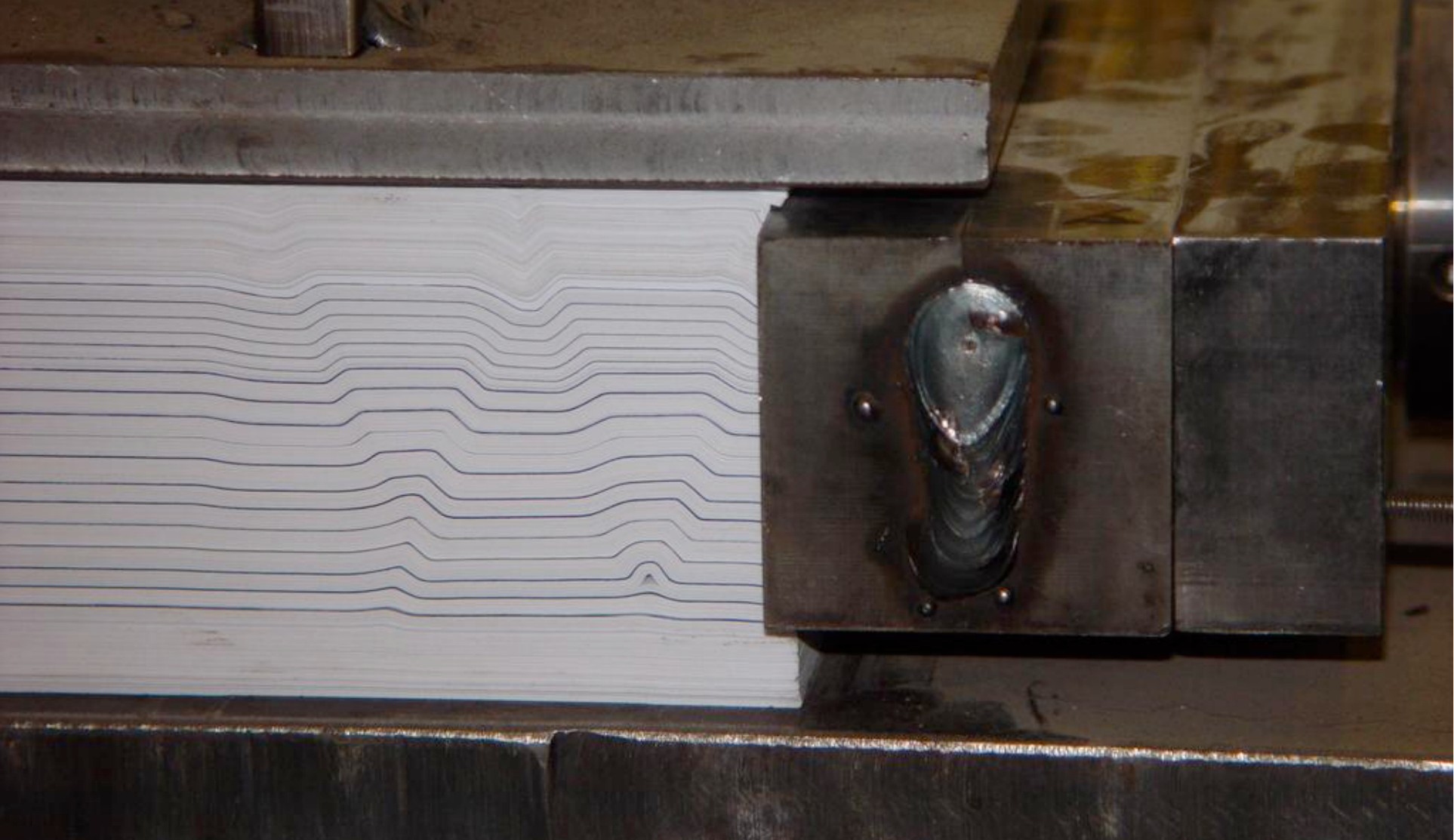}
\caption{An experiment similar to that of Fig.~\ref{fig:c5}, in which the foam is replaced by further sheets of paper, creating a much stiffer confinement. As a result the kink bands are shorter and more orthogonal to the layering. Experiment performed by Ahmer Wadee at Imperial College, London~\cite{WadeeHuntPeletier04}.}
\label{fig:Ahmer}
\end{figure}

This brought us to 
\begin{quote}
\textbf{Math question \#2: } How do kink bands form, and why?
\end{quote}
More specifically: what determines the geometry of the kink bands, such as the width of the kink band and its angle with respect to the layers? In the experiment of Fig.~\ref{fig:c5} the kink bands appear to arise as a sharpening of the global buckle that formed initially, and the size of the kink bands therefore is set by the initial buckle. That case is therefore not very interesting. In other experiments, such as in Fig.~\ref{fig:Ahmer}, the kink band width and angle seem to be independent of the size of the specimen. This is the case in which we want to understand how this geometry is determined.

\begin{figure}[htbp]
\centering
\includegraphics[height=6cm]{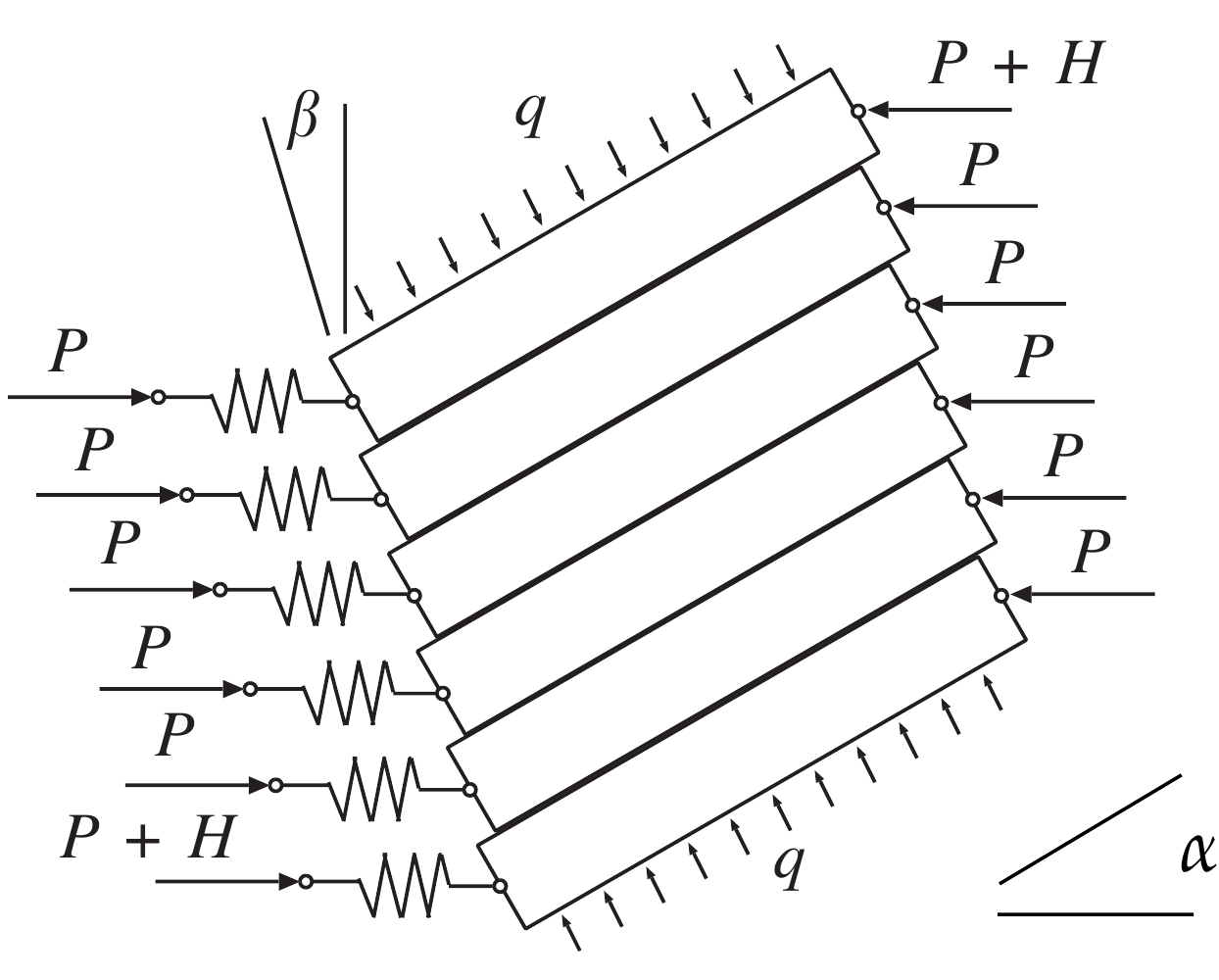}
\qquad\qquad\includegraphics[height=5.5cm]{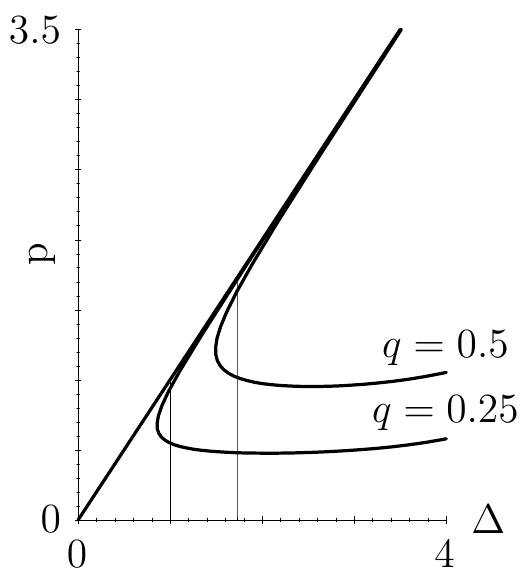}
\caption{Left: The setup of the simple kink band model. Right: Plot load $P$ versus shortening $\Delta$ for stationary points of the stack, for two values of overburden pressure $q$. The two vertical lines correspond to \emph{Maxwell displacements} (see text). Here $\mu=0.57$~\cite{HuntPeletierWadee00}.}
\label{fig:kinkbandsetup}
\end{figure}

Together with Giles Hunt and Ahmer Wadee we created a very simple mathematical model for kink bands (Rule \#4: \emph{Solve it, with the appliers}). The layers of paper in the kink band are abstracted into rigid blocks (see Fig.~\ref{fig:kinkbandsetup}), and the paper outside of the kink band is reduced to a single inline spring. The overburden pressure is represented by a simple lateral pressure $q$. The angle $\beta$ of the kink band is now fixed to the angle $\alpha$ of the layers by an assumption of non-expansion: the layers are assumed to preserve their thickness throughout deformation, and no opening between the layers is allowed. An important aspect of the model is \emph{friction} between the layers, characterized by the friction coefficient $\mu$, the ratio of the maximal friction force to the normal force.

Because of the friction, the system has an unusual structure. Many buckling problems have a finite \emph{buckling load}, at which the system becomes unstable. In this system however the friction causes the undeformed state to remain locally stable to arbitrarily high loads. As the load tends to infinity, the basin of stability shrinks to a point, and therefore the undeformed state becomes increasingly susceptible to small perturbations. The question becomes how to characterize this sensitivity to perturbations. 

We found a useful characterization in the concept of a \emph{Maxwell displacement}~\cite{HuntPeletierChampneysWoodsWadeeBuddLord00}, the smallest displacement $\Delta$ at which the energy of the undeformed state can be matched by the energy of a deformed state. The properties of this Maxwell displacement turn out to be favourable: the prediction of the kink band angle is reasonably accurate, and the prediction is remarkbly stable with respect to imperfections~\cite{HuntWadeePeletier99,HuntPeletierWadee00}. However, in this model the \emph{width} of the kink band is completely free. In order to fix that, the confinement, represented here by $q$, has to become stiffer with increasing compression~\cite{WadeeHuntPeletier04}.

\subsection{Revisting the sharp corners: voids}
\label{sec:voids}

We now return to the question of straight limbs and sharp corners. In such a sharp corner sometimes \emph{voids} arise; both Figs.~\ref{fig:c5} and~\ref{fig:Ahmer} show such voids between the layers. The simple pictorial argument described above of course does not allow any kind of prediction of the size of these voids, and therefore Tim Dodwell, a graduate student in Bath  with Giles Hunt and Chris Budd, studied a simplified problem aimed at characterizing the size of voids. 

\begin{figure}[htbp]
\centering
\includegraphics[height=8cm]{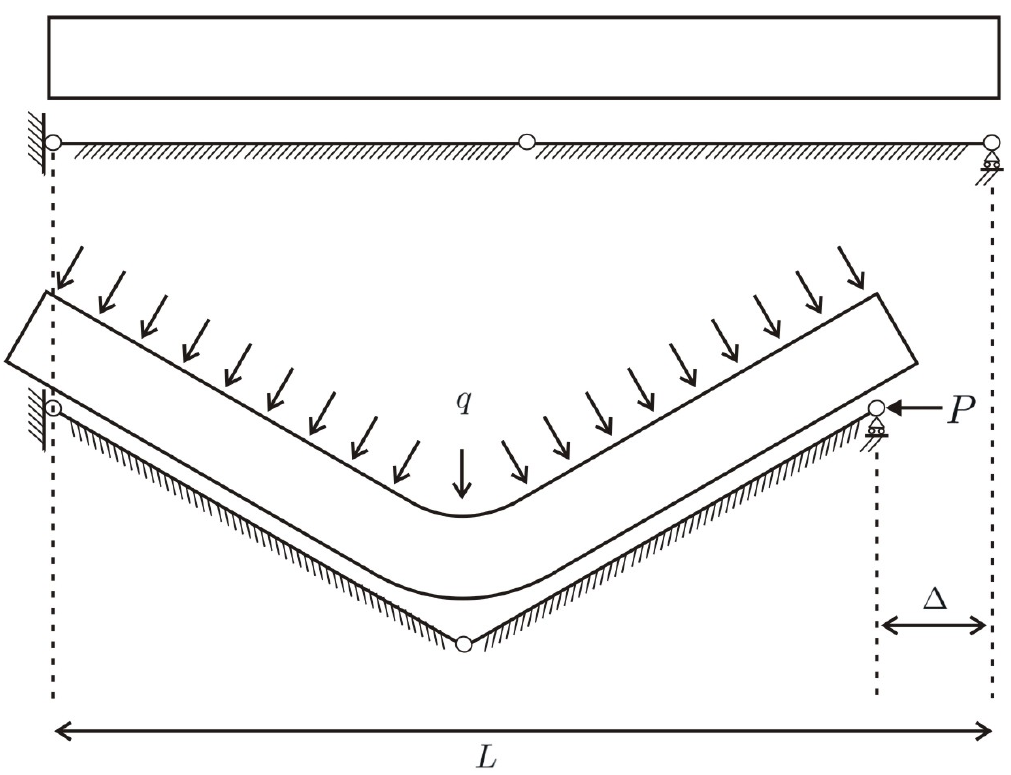}
\caption{Overburden pressure forces a single layer into a V-shaped obstacle~\cite{DodwellPeletierBuddHunt11}.}
\label{fig:voids}
\end{figure}

The setup is shown in Fig.~\ref{fig:voids}: an elastic layer is forced into a V-shaped obstacle by overburden pressure $q$ and possibly a lateral force $P$.

When $P=0$, this problem can be formulated as a constrained minimization problem, minimizing the functional
\[
V = \frac{B}{2}\int^{\infty}_{-\infty} \frac{w_{xx}^2}{(1 + w_x^2)^{5/2}} dx + q\int^{\infty}_{-\infty}(w-f)\,dx, \qquad \mbox{where } w\geq f. 
\]
Here $B$ is the bending stiffness of the rod, and $w:\R\to\R$ the vertical displacement; the obstacle is represented by the requirement $w\geq f$. Using ideas from convex function theory and constrained optimization, Tim proved a number of results:

\begin{theorem}[\cite{DodwellPeletierBuddHunt11}]
There exists a constrained global minimizer $w$, which is unique and convex. The void set $\{x:w(x)>f(x)\}$ is an interval, and its length scales as $(B/q)^{1/3}$. 
\end{theorem}

In addition, Tim derived the stationarity equation, characterized the contact forces, and gave a mechanical interpretation of the first integral of this equation. 

\medskip

This is a good start, and a marked progress with respect to Fig.~\ref{fig:sharp-corners}. The comparison with the experiment of Fig.~\ref{fig:c5}, however, is not altogether favourable. Indeed, that experiment shows an intriguing phenomenon: the voids do not happen every layer, but seem to be localized at black layers, and even at one in three or one in five black layers. 

The black layers are part of the experimental setup. In this experiment, we took a separate stack of paper and coloured the sides black with a felt pen. Then we introduced one of these coloured-side layers into the remaining stack, one every twenty-five layers of paper. 

So the position of the voids suggests that the process of introducing the black layers somehow caused the voids to happen there rather than somewhere else---perhaps as a consequence of reduced friction. But this does not explain why the voids happen not at every black layer, but on one in every five black layers---one in every 125 sheets of paper. 

\begin{figure}[htbp]
\centering
\includegraphics[height=5cm]{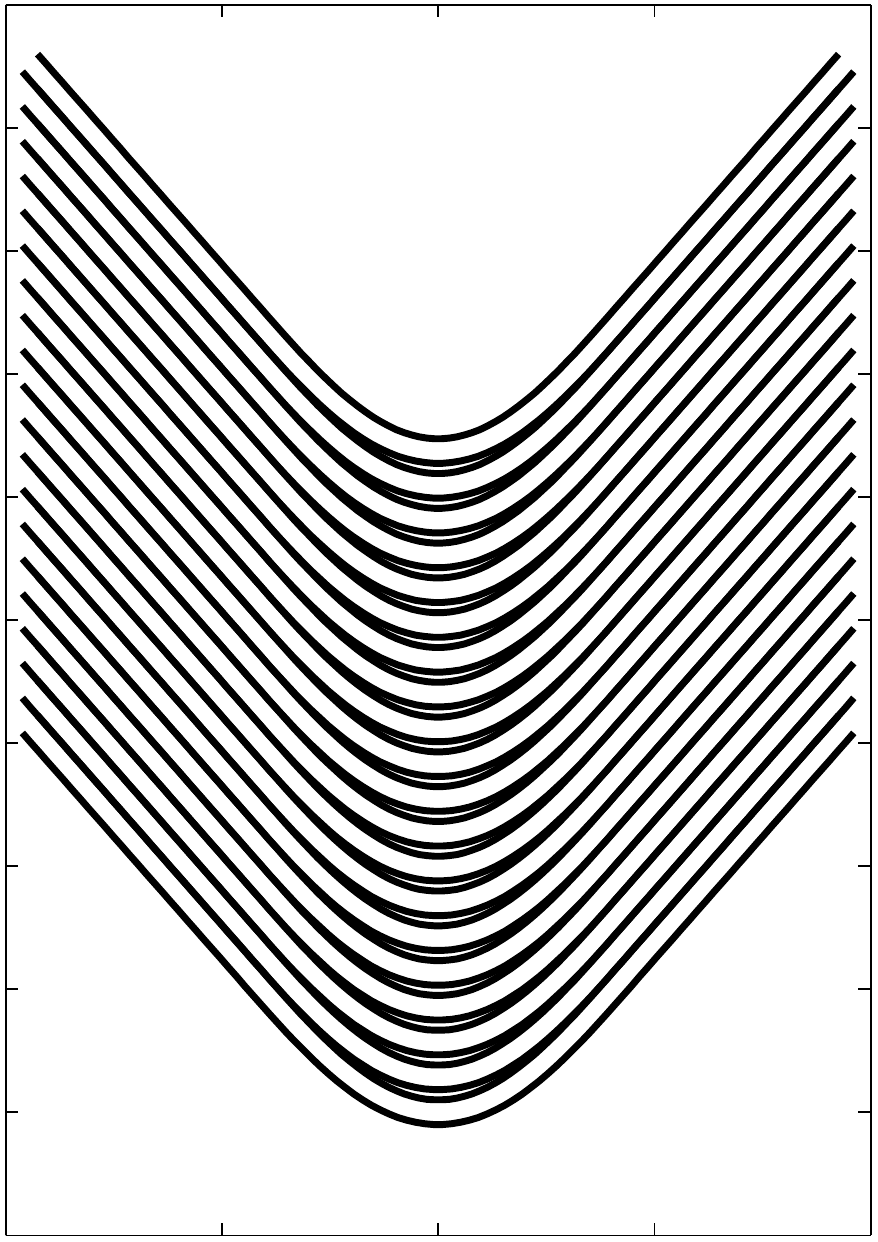}
\qquad
\includegraphics[height=5cm]{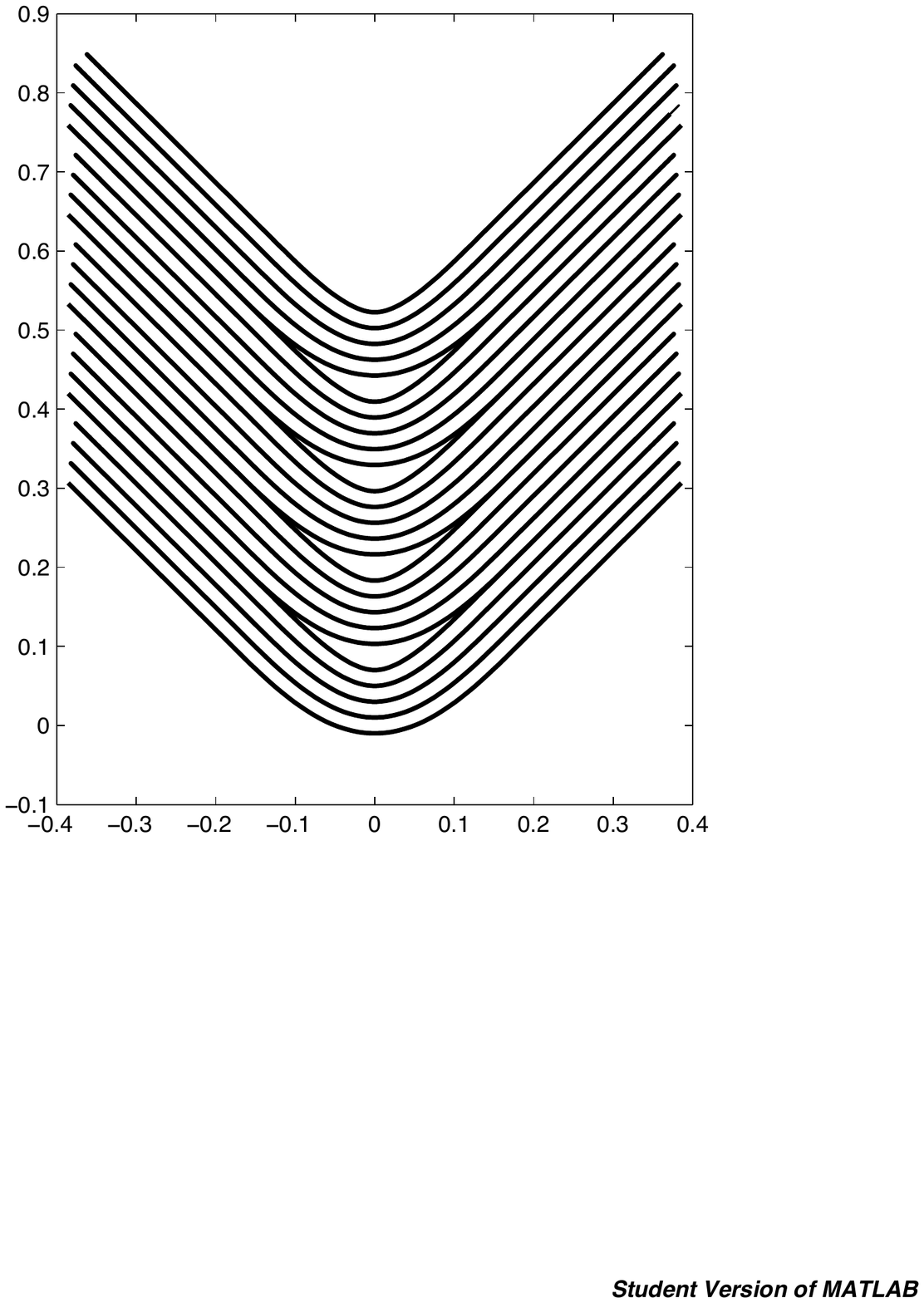}
\caption{In multilayer folding, voids might form between each two layers (left) or between each two packets of $n$ layers (right).}
\label{fig:multilayer}
\end{figure}

\begin{figure}[htb]
\centering
\includegraphics[width=10cm]{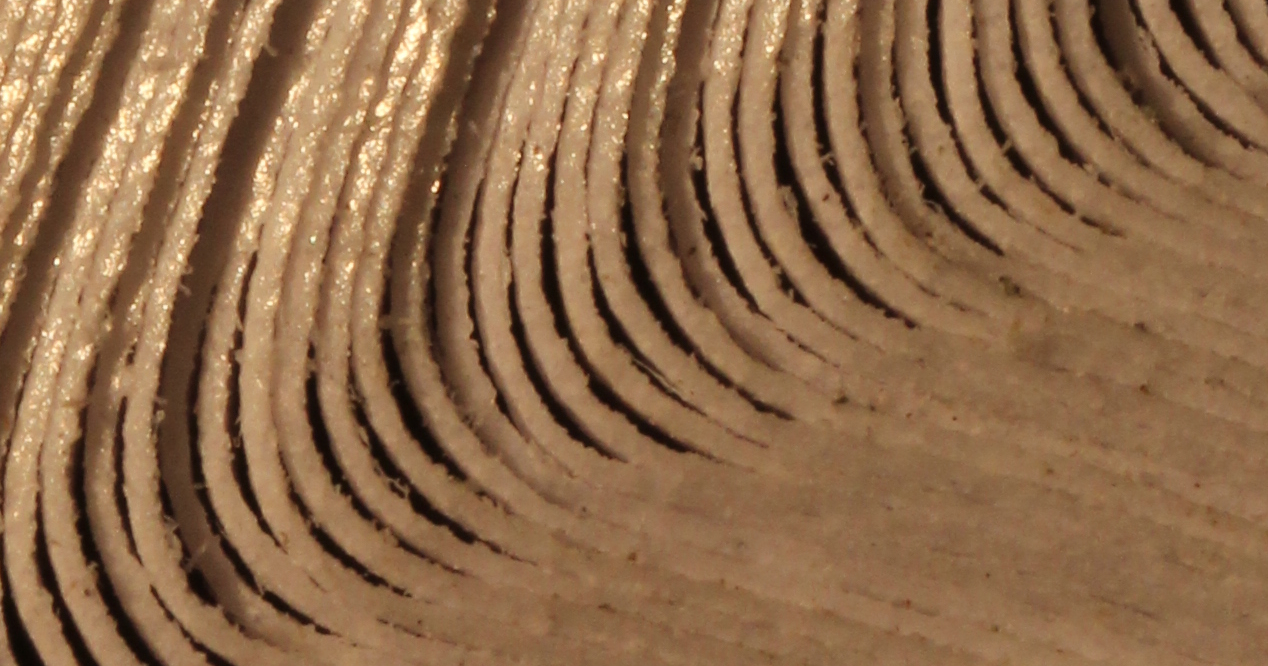}
\caption{In this experiment there appears to be a preference for void formation between each pair of layers. The photograph is highly enlarged, and shows individual layers of paper. Experiment performed by Tim Dodwell and Andrew Rhead at the University of Bath.}
\label{fig:expTim}
\end{figure}

Recently, Tim has been extending this model to cover multiple layers, again forced into a V-shaped obstacle, that might choose to deform separately (Fig~\ref{fig:multilayer}, left) or collectively (right). While it is too early to make strong predictions, there is a back-of-the-envelope calculation that suggests that for certain parameters layers might indeed choose to deform collectively in a pack of $n$ layers, with voids between two such packs, with an optimal $n$ that is larger than $1$ but smaller than $+\infty$~\cite{DodwellHuntPeletierBudd11TR}. Although that may be the case, other experiments done by Tim show preference for void formation between each pair of adjacent layers (Fig.~\ref{fig:expTim}). Obviously our understanding is still incomplete.

On a philosophical note, it is interesting to observe that while the paper experiments were set up to give us insight into rock folding, we are now trying to understand the behaviour of the paper experiments themselves. Hopefully this will eventually bring us back to the rocks again.

\subsection{Future work}

At this moment, there is at least one challenge that I would much like to address:
\begin{quote}
\textbf{Math question \#3: } characterize the layers of paper or rock \emph{as a material}.
\end{quote}
What I mean by this is the following. For homogeneous and isotropic elastic, viscous, and visco-elastic materials there is a well-defined constitutive theory, which provides the essential relationships between deformation and loading of those materials. For materials with some microstructure that remains within the same class the theory of homogenization provides rigorous derivations of material properties in the limit of large scale separation (see e.g.~\cite{Hornung97,CioranescuDonato99}). The `material' that we are considering here, however, is mixed in nature: the layers are elastic, while the interlayer mechanics introduces \emph{friction} into the picture. Such a combination has not been studied theoretically (although some related problems have, such as the case of \emph{rigid} `elasticity' with friction~\cite{ContiTheil05}). Other related problems are those of thin elastic bodies, for which rigorous derivations are available for a number of different energy regimes~\cite{FrieseckeJamesMuller06}. The case of elastic layers with friction remains open, however.

\section{Conclusion}

Looking back at the whole experience we can formulate a number of conclusions:
\begin{enumerate}
\item Suitable applied collaborators are very special, and very rare. One should cherish them!

\item The interaction with a particular applied problem does not necessarily result in an interesting mathematical problem \emph{right away}. The first two results described above were mathematically uninteresting: the pictorial argument isn't even  mathematical, and the kind-band model in the form of a stack borders on the trivial, from a mathematical point of view. However, as time progressed, we were able to move to more interesting mathematical models, such as the variational model in Section~\ref{sec:voids}, as answers to the applied questions. 

This phenomenon is quite general, and implies a lesson learnt: finding good mathematical problems takes time, energy, and perseverance---but the result is worth it.

\item This experience is an example of a collaboration in which both the application side and the mathematical side took home good and interesting results. I believe that this is often possible, and I also believe that one should aim for such a collaboration. This will take time and energy---see above---but it is worth it.
\end{enumerate}

In retrospect I realize that during my PhD period in Delft and Amsterdam, Hans had a very similar collaboration going with a group of petroleum engineers in Delft, notably Hans Bruining, and a group of subsurface soil engineers, including Sjoerd van de Zee and Pieter Raats. As I remarked in the beginning, we never really explicitly discussed \emph{how} one handles such a collaboration. But I have been very fortunate to have seen this way of working in action, and my subsequent work in applied mathematics---in mathematical geology and other areas---have been strongly influenced by it. For me, Hans has been a wonderful teacher.

\bibliography{ref}
\bibliographystyle{plain}

\end{document}